# Finite Mixtures of Multivariate Skew Laplace Distributions

Fatma Zehra Doğru[1]*, Y. Murat Bulut[2] and Olcay Arslan[3]

[1] Giresun University, Faculty of Arts and Sciences, Department of Statistics, 28100 Giresun/Turkey.
(E-mail: fatma.dogru@giresun.edu.tr)
[2] Eskisehir Osmangazi University, Faculty of Science and Letters, Department of Statistics, 26480 Eskisehir/Turkey.(E-mail: ymbulut@ogu.edu.tr)
[3] Ankara University, Faculty of Science, Department of Statistics, 06100 Ankara/Turkey.
(E-mail: oarslan@ankara.edu.tr)

**Abstract**

In this paper, we propose finite mixtures of multivariate skew Laplace distributions to model both skewness and heavy-tailedness in the heterogeneous data sets. The maximum likelihood estimators for the parameters of interest are obtained by using the EM algorithm. We give a small simulation study and a real data example to illustrate the performance of the proposed mixture model.

## 1. Introduction

Finite mixture models are used for modeling heterogeneous data sets as a result of their flexibility. These models are commonly applied in fields such as classification, cluster and latent class analysis, density estimation, data mining, image analysis, genetics, medicine, pattern recognition etc. (see for more detail, Titterington et al. (1985), McLachlan and Basford (1988), McLachlan and Peel (2000), Bishop (2006), Frühwirth-Schnatter (2006)).

Generally, due to its tractability and wide applicability, the distribution of mixture model components is assumed to be normal. However, in practice, the data sets may be asymmetric and/or heavy-tailed. For this, there are some studies in literature for multivariate mixture modeling using the asymmetric and/or heavy-tailed distributions. Some of these studies can be summarized as follows. Peel and McLachlan (2000) proposed finite mixtures of multivariate t distributions as a robust extension of the multivariate normal mixture model (McLachlan and Basford (1988)), Lin (2009) introduced multivariate skew normal mixture models, Pyne et al. (2009) and Lin (2010) proposed finite mixtures of the restricted and unrestricted variants of multivariate skew t distributions of Sahu et al. (2003), Cabral et al. (2012) proposed multivariate mixture modeling based on the skew-normal independent distributions and Lin et al. (2014) introduced a flexible mixture modeling based on the skew-t-normal distribution.

In the multivariate analysis, the multivariate skew normal (MSN) (Azzalini and Dalla Valle (1996), Gupta et al. (2004) and Arellano-Valle and Genton (2005)) distribution has been proposed as an alternative to the multivariate normal (MN) distribution to deal with skewness in the data. However, since MSN distribution is not heavy-tailed some alternative heavy-tailed skew distributions are needed to model skewness and heavy-tailedness. One of the examples of heavy-tailed skew distribution is the multivariate skew t (MST) distribution defined by Azzalini and Capitanio (2003) and Gupta (2003). The other heavy-tailed skew distribution is the multivariate skew Laplace (MSL) distribution proposed by Arslan (2010). The advantage of this distribution is that the MSL distribution has the less number of parameters than MST distribution and it has the same number of parameters with the MSN distribution. Concerning the finite mixtures of distributions, since the mixtures of MN distributions are not able to model skewness, finite mixtures of MSN distributions have been proposed by Lin (2009) to model such data sets. In this study, we explore the finite mixtures of MSL distributions as an alternative to the finite mixtures of MSN distributions to deal with both skewness and heavy-tailedness in the heterogeneous data sets.

The rest of this paper is organized as follows. In Section 2, we briefly summarize some properties of MSL distribution. In Section 3, we present the mixtures of MSL distributions and give the Expectation-Maximization (EM) algorithm to obtain the maximum likelihood (ML) estimators of proposed mixture model. In Section 4, we give the empirical information matrix of MSL distribution to compute the standard errors of proposed estimators. In the Application Section, we give a small simulation study and a real data example to illustrate the performance of proposed mixture model. Some conclusions are given in Section 6.

## 2. Multivariate skew Laplace distribution

A $p$-dimensional random vector $\boldsymbol{Y} \in R^p$ is said to have MSL distribution ($Y \sim MSL_p(\boldsymbol{\mu}, \Sigma, \boldsymbol{\gamma})$) which is given by Arslan (2010) if it has the following probability density function (pdf)

$$f_{MSL}(\boldsymbol{y}; \boldsymbol{\mu}, \Sigma, \boldsymbol{\gamma}) = \frac{|\Sigma|^{-\frac{1}{2}}}{2^p \pi^{\frac{p-1}{2}} \alpha \Gamma\left(\frac{p+1}{2}\right)} \exp\left\{-\alpha\sqrt{(\boldsymbol{y}-\boldsymbol{\mu})^T\Sigma^{-1}(\boldsymbol{y}-\boldsymbol{\mu})} + (\boldsymbol{y}-\boldsymbol{\mu})^T\Sigma^{-1}\boldsymbol{\gamma}\right\}, \quad (1)$$

where $\alpha = \sqrt{1 + \boldsymbol{\gamma}^T\Sigma^{-1}\boldsymbol{\gamma}}$, $\boldsymbol{\mu} \in R^p$ is a location parameter, $\boldsymbol{\gamma} \in R^p$ is a skewness parameter and $\Sigma$ is a positive definite scatter matrix.

**Proposition 1.** The characteristic function of $MSL_n(\boldsymbol{\mu}, \Sigma, \boldsymbol{\gamma})$ is

$$\Phi_Y(t) = e^{it^T\mu}[1 + t^T\Sigma t - 2it^T\boldsymbol{\gamma}]^{-(p+1)/2}, \quad t\epsilon R^p.$$

**Proof.** Since we know conditional distribution of $\boldsymbol{Y}$ given V, we get characteristic function as follows

$$\begin{aligned}\Phi_Y(t) &= E\left(E_{Y|V}\left(e^{it^TY}\right)\right) = e^{it^T\mu}E_V\left(e^{-V^{-1}\left(\frac{t^T\Sigma t}{2} - it^T\gamma\right)}\right) \\ &= e^{it^T\mu}\int_0^\infty e^{v^{-1}\left(\frac{t^T\Sigma t}{2} - it^T\gamma\right)}g(v)dv \\ &= e^{it^T\mu}[1 + t^T\Sigma t - 2it^T\boldsymbol{\gamma}]^{-(n+1)/2}.\end{aligned}$$ ∎

If $\boldsymbol{Y} \sim MSL_p(\mu, \Sigma, \gamma)$ then the expectation and variance of $\boldsymbol{Y}$ is

$$E(\boldsymbol{Y}) = \boldsymbol{\mu} + (p+1)\boldsymbol{\gamma},$$
$$Var(\boldsymbol{Y}) = (p+1)(\Sigma + 2\boldsymbol{\gamma}\boldsymbol{\gamma}^T).$$

The MSL distribution can be obtained as a variance-mean mixture of MN distribution and inverse gamma (IG) distribution. The variance-mean mixture representation is given as follows

$$\boldsymbol{Y} = \boldsymbol{\mu} + V^{-1}\boldsymbol{\gamma} + \sqrt{V^{-1}}\Sigma^{1/2}\boldsymbol{X} \qquad (2)$$

where $X \sim N_p(\mathbf{0}, I_p)$ and $V \sim IG(\frac{p+1}{2}, \frac{1}{2})$. Note that if $\boldsymbol{\gamma} = \mathbf{0}$, the density function of $Y$ reduces to the density function of symmetric multivariate Laplace distribution given by Naik and Plungpongpun (2006). Also, the conditional distribution of $Y$ given $V = v$ will be

$$Y|v \sim N(\boldsymbol{\mu} + v^{-1}\boldsymbol{\gamma}, v^{-1}\Sigma).$$

The joint density function of $Y$ and $V$ is

$$f(\mathbf{y}, v) = \frac{|\Sigma|^{-\frac{1}{2}} e^{(\mathbf{y}-\boldsymbol{\mu})^T \Sigma^{-1} \boldsymbol{\gamma}}}{2^p \pi^{\frac{p-1}{2}} \alpha \Gamma\left(\frac{p+1}{2}\right)} \left\{ v^{-3/2} e^{-\frac{1}{2}\{(\mathbf{y}-\boldsymbol{\mu})^T \Sigma^{-1}(\mathbf{y}-\boldsymbol{\mu})v + (1+\boldsymbol{\gamma}^T \Sigma^{-1} \boldsymbol{\gamma})v^{-1}\}} \right\}.$$

Then, we have the following conditional density function of $V$ given $Y$

$$f(v|\mathbf{y}) = \frac{\alpha}{\sqrt{2\pi}} e^{\alpha \sqrt{(\mathbf{y}-\boldsymbol{\mu})^T \Sigma^{-1}(\mathbf{y}-\boldsymbol{\mu})}} v^{-3/2} e^{-\frac{1}{2}\{(\mathbf{y}-\boldsymbol{\mu})^T \Sigma^{-1}(\mathbf{y}-\boldsymbol{\mu})v + \alpha^2 v^{-1}\}}, \quad v > 0. \quad (3)$$

**Proposition 2.** Using the conditional density function given in (3), the conditional expectations can be obtained as follows

$$E(V|\mathbf{y}) = \frac{\sqrt{1 + \boldsymbol{\gamma}^T \Sigma^{-1} \boldsymbol{\gamma}}}{\sqrt{(\mathbf{y} - \boldsymbol{\mu})^T \Sigma^{-1}(\mathbf{y} - \boldsymbol{\mu})}}, \quad (4)$$

$$E(V^{-1}|\mathbf{y}) = \frac{1 + \sqrt{(1 + \boldsymbol{\gamma}^T \Sigma^{-1} \boldsymbol{\gamma})(\mathbf{y} - \boldsymbol{\mu})^T \Sigma^{-1}(\mathbf{y} - \boldsymbol{\mu})}}{1 + \boldsymbol{\gamma}^T \Sigma^{-1} \boldsymbol{\gamma}}. \quad (5)$$

Note that these conditional expectations will be used in the EM algorithm given in Subsection 3.1.

## 3. Finite mixtures of the MSL distributions

Let $\mathbf{y}_1, \dots, \mathbf{y}_n$ be a $p$-dimensional random sample which come from a $g$-component mixtures of MSL distributions. The pdf of a $g$-component finite mixtures of MSL distributions is given by

$$f(\mathbf{y}|\boldsymbol{\Theta}) = \sum_{i=1}^{g} \pi_i f(\mathbf{y}; \boldsymbol{\mu}_i, \Sigma_i, \boldsymbol{\gamma}_i), \quad (6)$$

where $\pi_i$ denotes the mixing probability with $\sum_{i=1}^{g} \pi_i = 1$, $0 \leq \pi_i \leq 1$, $f(\mathbf{y}; \boldsymbol{\mu}_i, \Sigma_i, \boldsymbol{\lambda}_i)$ represents the pdf of the $ith$ component (pdf of the MSL distribution) given in (1) and $\boldsymbol{\Theta} = (\pi_1, \dots, \pi_g, \boldsymbol{\mu}_1, \dots, \boldsymbol{\mu}_g, \Sigma_1, \dots, \Sigma_g, \boldsymbol{\gamma}_1, \dots, \boldsymbol{\gamma}_g)^T$ is the unknown parameter vector.

### 3.1 ML estimation

The ML estimator of $\boldsymbol{\Theta}$ can be found by maximizing the following log-likelihood function

$$\ell(\boldsymbol{\Theta}) = \sum_{j=1}^{n} \log \left( \sum_{i=1}^{g} \pi_i f(\mathbf{y}_j; \boldsymbol{\mu}_i, \Sigma_i, \boldsymbol{\gamma}_i) \right). \quad (7)$$

However, there is not an explicit maximizer of (7). Therefore, in general, the EM algorithm (Dempster et al. (1977)) is used to obtain the ML estimator of $\boldsymbol{\Theta}$. Here, we will use the following EM algorithm.

Let $\boldsymbol{z}_j = (z_{1j}, \ldots, z_{gj})^T$ be the latent variables with

$$z_{ij} = \begin{cases} 1, & \text{if } j^{th} \text{ observation belongs to } i^{th} \text{ component} \\ 0, & \text{otherwise} \end{cases} \quad (8)$$

where $j = 1, \ldots, n$ and $i = 1, \ldots, g$. To implement the steps of the EM algorithm, we will use the stochastic representation of the MSL distribution given in (2). Then, the hierarchical representation for the mixtures of MSL distributions will be

$$\boldsymbol{Y}_j | v_j, z_{ij} = 1 \sim N(\boldsymbol{\mu} + v_j^{-1} \boldsymbol{\gamma}, v_j^{-1} \Sigma),$$
$$V_j | z_{ij} = 1 \sim IG\left(\frac{p+1}{2}, \frac{1}{2}\right). \quad (9)$$

Let $(\boldsymbol{y}, \boldsymbol{v}, \boldsymbol{z})$ be the complete data, where $\boldsymbol{y} = (\boldsymbol{y}_1^T, \ldots, \boldsymbol{y}_n^T)^T$, $\boldsymbol{v} = (v_1, \ldots, v_n)$ and $\boldsymbol{z} = (\boldsymbol{z}_1, \ldots, \boldsymbol{z}_n)^T$. Using the hierarchical representation given above and ignoring the constants, the complete data log-likelihood function can be written by

$$\ell_c(\boldsymbol{\Theta}; \boldsymbol{y}, \boldsymbol{v}, \boldsymbol{z}) = \sum_{j=1}^n \sum_{i=1}^n z_{ij} \left\{ \log \pi_i - \frac{1}{2} \log |\Sigma_i| + (\boldsymbol{y}_j - \boldsymbol{\mu}_i)^T \Sigma_i^{-1} \boldsymbol{\gamma}_i \right.$$
$$\left. - \frac{1}{2} v_j (\boldsymbol{y}_j - \boldsymbol{\mu}_i)^T \Sigma_i^{-1} (\boldsymbol{y}_j - \boldsymbol{\mu}_i) - \frac{1}{2} \boldsymbol{\gamma}_i^T \Sigma_i^{-1} \boldsymbol{\gamma}_i v_j^{-1} - \frac{1}{2} (3 \log v_j + v_j^{-1}) \right\}. \quad (10)$$

To overcome the latency of the latent variables given in (10), we have to take the conditional expectation of the complete data log-likelihood function given the observed data $\boldsymbol{y}_j$

$$E(\ell_c(\boldsymbol{\Theta}; \boldsymbol{y}, \boldsymbol{v}, \boldsymbol{z}) | \boldsymbol{y}_j) = \sum_{j=1}^n \sum_{i=1}^g E(z_{ij} | \boldsymbol{y}_j) \left\{ \log \pi_i - \frac{1}{2} \log |\Sigma_i| - (\boldsymbol{y}_j - \boldsymbol{\mu}_i)^T \Sigma_i^{-1} \boldsymbol{\gamma}_i \right.$$
$$\left. - \frac{1}{2} E(V_j | \boldsymbol{y}_j) (\boldsymbol{y}_j - \boldsymbol{\mu}_i)^T \Sigma_i^{-1} (\boldsymbol{y}_j - \boldsymbol{\mu}_i) - \frac{1}{2} \boldsymbol{\gamma}_i^T \Sigma_i^{-1} \boldsymbol{\gamma}_i E(V_j^{-1} | \boldsymbol{y}_j) \right\}. \quad (11)$$

Since the last part of the complete data log-likelihood function does not include the parameters, the above conditional expectation of the complete data log-likelihood function only consists of the necessary conditional expectations.

Note that the conditional expectations $E(V_j | \boldsymbol{y}_j)$ and $E(V_j^{-1} | \boldsymbol{y}_j)$ can be calculated using the conditional expectations given in (4) and (5) and the conditional expectation $E(z_{ij} | \boldsymbol{y}_j)$ can be computed using the classical theory of mixture modeling.

**EM algorithm:**

**1.** Set initial parameter estimate $\boldsymbol{\Theta}^{(0)}$ and a stopping rule $\Delta$.
**2. E-Step:** Compute the following conditional expectations for $k = 0,1,2,\ldots$ iteration

$$\hat{z}_{ij}^{(k)} = E(z_{ij} | \boldsymbol{y}_j, \widehat{\boldsymbol{\Theta}}^{(k)}) = \frac{\hat{\pi}_i^{(k)} f\left(\boldsymbol{y}_j; \widehat{\boldsymbol{\mu}}_i^{(k)}, \widehat{\Sigma}_i^{(k)}, \widehat{\boldsymbol{\gamma}}_i^{(k)}\right)}{f(\boldsymbol{y}_j; \widehat{\boldsymbol{\Theta}}^{(k)})}, \quad (12)$$

$$\hat{v}_{1ij}^{(k)} = E(V_j|\mathbf{y}_j, \widehat{\boldsymbol{\Theta}}^{(k)}) = \frac{\sqrt{1 + \widehat{\boldsymbol{\gamma}}_i^{(k)T}\widehat{\Sigma}_i^{(k)-1}\widehat{\boldsymbol{\gamma}}_i^{(k)}}}{\sqrt{\left(\mathbf{y}_j - \widehat{\boldsymbol{\mu}}_i^{(k)}\right)^T \widehat{\Sigma}_i^{(k)-1} \left(\mathbf{y}_j - \widehat{\boldsymbol{\mu}}_i^{(k)}\right)}}, \tag{13}$$

$$\hat{v}_{2ij}^{(k)} = E(V_j^{-1}|\mathbf{y}_j, \widehat{\boldsymbol{\Theta}}^{(k)})$$

$$= \frac{1 + \sqrt{\left(1 + \widehat{\boldsymbol{\gamma}}_i^{(k)T}\widehat{\Sigma}_i^{(k)-1}\widehat{\boldsymbol{\gamma}}_i^{(k)}\right)\left(\mathbf{y}_j - \widehat{\boldsymbol{\mu}}_i^{(k)}\right)^T \widehat{\Sigma}_i^{(k)-1}\left(\mathbf{y}_j - \widehat{\boldsymbol{\mu}}_i^{(k)}\right)}}{1 + \widehat{\boldsymbol{\gamma}}_i^{(k)T}\widehat{\Sigma}_i^{(k)-1}\widehat{\boldsymbol{\gamma}}_i^{(k)}}. \tag{14}$$

Then, we form the following objective function

$$Q(\boldsymbol{\Theta}; \widehat{\boldsymbol{\Theta}}^{(k)}) = \sum_{j=1}^n \sum_{i=1}^g \hat{z}_{ij}^{(k)} \left\{ \log \pi_i - \frac{1}{2}\log|\Sigma_i| - (\mathbf{y}_j - \boldsymbol{\mu}_i)^T \Sigma_i^{-1} \boldsymbol{\gamma}_i \right. $$
$$\left. - \frac{1}{2}\hat{v}_{1ij}^{(k)}(\mathbf{y}_j - \boldsymbol{\mu}_i)^T \Sigma_i^{-1}(\mathbf{y}_j - \boldsymbol{\mu}_i) - \frac{1}{2}\hat{v}_{2ij}^{(k)} \boldsymbol{\gamma}_i^T \Sigma_i^{-1} \boldsymbol{\gamma}_i \right\}. \tag{15}$$

**3. M-step:** Maximize the $Q(\boldsymbol{\Theta}; \widehat{\boldsymbol{\Theta}}^{(k)})$ with respect to $\boldsymbol{\Theta}$ to get the $(k+1)th$ parameter estimates for the parameters. This maximization yields the following updating equations

$$\hat{\pi}_i^{(k+1)} = \frac{\sum_{j=1}^n \hat{z}_{ij}^{(k)}}{n}, \tag{16}$$

$$\hat{\boldsymbol{\mu}}_i^{(k+1)} = \frac{\sum_{j=1}^n \hat{z}_{ij}^{(k)} \hat{v}_{1ij}^{(k)} \mathbf{y}_j - \sum_{j=1}^n \hat{z}_{ij}^{(k)} \widehat{\boldsymbol{\gamma}}_i^{(k)}}{\sum_{j=1}^n \hat{z}_{ij}^{(k)} \hat{v}_{1ij}^{(k)}}, \tag{17}$$

$$\widehat{\boldsymbol{\gamma}}_i^{(k+1)} = \frac{\left(\sum_{j=1}^n \hat{z}_{ij}^{(k)} \hat{v}_{1ij}^{(k)}\right)\left(\sum_{j=1}^n \hat{z}_{ij}^{(k)} \mathbf{y}_j\right) - \left(\sum_{j=1}^n \hat{z}_{ij}^{(k)}\right)\left(\sum_{j=1}^n \hat{z}_{ij}^{(k)} \hat{v}_{1ij}^{(k)} \mathbf{y}_j\right)}{\left(\sum_{j=1}^n \hat{z}_{ij}^{(k)} \hat{v}_{1ij}^{(k)}\right)\left(\sum_{j=1}^n \hat{z}_{ij}^{(k)} \hat{v}_{2ij}^{(k)}\right) - \left(\sum_{j=1}^n \hat{z}_{ij}^{(k)}\right)^2}, \tag{18}$$

$$\widehat{\Sigma}_i^{(k+1)} = \frac{\sum_{j=1}^n \hat{z}_{ij}^{(k)} \hat{v}_{1ij}^{(k)} \left(\mathbf{y}_j - \widehat{\boldsymbol{\mu}}_i^{(k)}\right)\left(\mathbf{y}_j - \widehat{\boldsymbol{\mu}}_i^{(k)}\right)^T - \widehat{\boldsymbol{\gamma}}_i^{(k)} \widehat{\boldsymbol{\gamma}}_i^{(k)T} \sum_{j=1}^n \hat{z}_{ij}^{(k)} \hat{v}_{2ij}^{(k)}}{\sum_{j=1}^n \hat{z}_{ij}^{(k)}}. \tag{19}$$

**4.** Repeat E and M steps until the convergence rule $\|\widehat{\boldsymbol{\Theta}}^{(k+1)} - \widehat{\boldsymbol{\Theta}}^{(k)}\| < \Delta$ is obtained. Alternatively, the absolute difference of the actual log-likelihood $\|\ell(\widehat{\boldsymbol{\Theta}}^{(k+1)}) - \ell(\widehat{\boldsymbol{\Theta}}^{(k)})\| < \Delta$ or $\|\ell(\widehat{\boldsymbol{\Theta}}^{(k+1)})/\ell(\widehat{\boldsymbol{\Theta}}^{(k)}) - 1\| < \Delta$ can be used as a stopping rule (see Dias and Wedel (2004)).

### 3.2 Initial values

To determine the initial values for the EM algorithm, we will use the same procedure given in Lin (2009). The steps of the selecting initial values are given as follows.

i) Perform the K-means clustering algorithm (Hartigan and Wong (1979)).
ii) Initialize the component labels $\hat{\mathbf{z}}_j^{(0)} = \{z_{ij}\}_{i=1}^g$ according to the K-means clustering results.
iii) The initial values of mixing probabilities, component locations and component scale variances can be set as

$$\hat{\pi}_i^{(0)} = \frac{\sum_{j=1}^n \hat{z}_{ij}^{(0)}}{n}, \quad \hat{\boldsymbol{\mu}}_i^{(0)} = \frac{\sum_{j=1}^n \hat{z}_{ij}^{(0)} \boldsymbol{y}_j}{\sum_{j=1}^n \hat{z}_{ij}^{(0)}},$$

$$\hat{\Sigma}_i^{(0)} = \frac{\sum_{j=1}^n \hat{z}_{ij}^{(0)} \left(\boldsymbol{y}_j - \hat{\boldsymbol{\mu}}_i^{(0)}\right)\left(\boldsymbol{y}_j - \hat{\boldsymbol{\mu}}_i^{(0)}\right)^T}{\sum_{j=1}^n \hat{z}_{ij}^{(0)}}.$$

iv) For the skewness parameters, use the skewness coefficient vector of each clusters.

## 4. The empirical information matrix

We will compute the standard errors of ML estimators using the information based method given by Basford et al. (1997). Here, we will use the inverse of the empirical information matrix to have an approximation to the asymptotic covariance matrix of estimators. The information matrix is

$$\hat{I}_e = \sum_{j=1}^n \hat{\boldsymbol{s}}_j \hat{\boldsymbol{s}}_j^T, \tag{20}$$

where $\hat{\boldsymbol{s}}_j = E_{\hat{\Theta}}\left(\frac{\partial \ell_{cj}(\Theta; \boldsymbol{y}_j, \boldsymbol{v}_j, \boldsymbol{z}_j)}{\partial \Theta}\bigg| \boldsymbol{y}_j\right), j = 1, \ldots, n$ are the individual scores and $\ell_{cj}(\Theta; \boldsymbol{y}_j, \boldsymbol{v}_j, \boldsymbol{z}_j)$ is the complete data log-likelihood function for the $jth$ observation. Then, the elements of the score vector $\hat{\boldsymbol{s}}_j$ is $\left(\hat{\boldsymbol{s}}_{j,\pi_1}, \ldots, \hat{\boldsymbol{s}}_{j,\pi_{g-1}}, \hat{\boldsymbol{s}}_{j,\boldsymbol{\mu}_1}, \ldots, \hat{\boldsymbol{s}}_{j,\boldsymbol{\mu}_g}, \hat{\boldsymbol{s}}_{j,\sigma_1}, \ldots, \hat{\boldsymbol{s}}_{j,\sigma_g}, \hat{\boldsymbol{s}}_{j,\boldsymbol{\gamma}_1}, \ldots, \hat{\boldsymbol{s}}_{j,\boldsymbol{\gamma}_g}\right)^T$. After some straightforward algebra, we obtain the following equations

$$\hat{\boldsymbol{s}}_{j,\pi_r} = \frac{\hat{z}_{rj}}{\hat{\pi}_r} - \frac{\hat{z}_{gj}}{\hat{\pi}_g}, \qquad r = 1, \ldots, g-1, \tag{21}$$

$$\hat{\boldsymbol{s}}_{j,\boldsymbol{\mu}_i} = \hat{z}_{ij} \hat{\Sigma}_i^{-1}\left(\hat{v}_{1ij}(\boldsymbol{y}_j - \hat{\boldsymbol{\mu}}_i) - \hat{\boldsymbol{\gamma}}_i\right), \tag{22}$$

$$\hat{\boldsymbol{s}}_{j,\sigma_i} = vech\left(\hat{z}_{ij}\left\{-\left(\hat{\Sigma}_i^{-1} - \hat{v}_{1ij}\hat{\Sigma}_i^{-1}(\boldsymbol{y}_j - \hat{\boldsymbol{\mu}}_i)(\boldsymbol{y}_j - \hat{\boldsymbol{\mu}}_i)'\hat{\Sigma}_i^{-1} - \hat{v}_{2ij}\hat{\Sigma}_i^{-1}\hat{\boldsymbol{\gamma}}_i\hat{\boldsymbol{\gamma}}_i^T\hat{\Sigma}_i^{-1}\right)\right.\right.$$
$$\left.\left. + \frac{1}{2} diag\left(\hat{\Sigma}_i^{-1} - \hat{v}_{1ij}\hat{\Sigma}_i^{-1}(\boldsymbol{y}_j - \hat{\boldsymbol{\mu}}_i)(\boldsymbol{y}_j - \hat{\boldsymbol{\mu}}_i)'\hat{\Sigma}_i^{-1}\right.\right.\right. \tag{23}$$
$$\left.\left.\left. - \hat{v}_{2ij}\hat{\Sigma}_i^{-1}\hat{\boldsymbol{\gamma}}_i\hat{\boldsymbol{\gamma}}_i^T\right)\right\}\right),$$

$$\hat{\boldsymbol{s}}_{j,\boldsymbol{\gamma}_i} = \hat{z}_{ij}\hat{\Sigma}_i^{-1}\left((\boldsymbol{y}_j - \hat{\boldsymbol{\mu}}_i) - \hat{v}_{2ij}\hat{\boldsymbol{\gamma}}_i\right). \tag{24}$$

Thus, the standard errors of $\hat{\Theta}$ will be found by using the square root of the diagonal elements of the inverse of (20).

## 5. Applications

In this section, we will illustrate the performance of proposed mixture model with a small simulation study and a real data example. All computations for simulation study and real data example were conducted using MATLAB R2013a. For all computations, the stopping rule $\Delta$ is taken as $10^{-6}$.

### 5.1 Simulation study

In the simulation study, the data are generated from the following two-component mixtures of MSL distributions

$$f(\mathbf{y}_j|\Theta) = \pi_1 f_p(\mathbf{y}_j; \boldsymbol{\mu}_1, \Sigma_1, \boldsymbol{\gamma}_1) + (1-\pi_1) f_p(\mathbf{y}_j; \boldsymbol{\mu}_2, \Sigma_2, \boldsymbol{\gamma}_2),$$

where

$$\boldsymbol{\mu}_i = (\mu_{i1}, \mu_{i2})^T, \quad \Sigma_i = \begin{bmatrix} \sigma_{i,11} & \sigma_{i,12} \\ \sigma_{i,21} & \sigma_{i,22} \end{bmatrix}, \quad \boldsymbol{\gamma}_i = (\gamma_{i1}, \gamma_{i2})^T, i = 1,2$$

with the parameter values

$$\boldsymbol{\mu}_1 = (2,2)^T, \quad \boldsymbol{\mu}_2 = (-2,-2)^T, \quad \Sigma_1 = \Sigma_2 = \begin{bmatrix} 1.5 & 0 \\ 0 & 1.5 \end{bmatrix},$$
$$\boldsymbol{\gamma}_1 = (1,1)^T, \quad \boldsymbol{\gamma}_2 = (-1,-1)^T, \quad \pi_1 = 0.6.$$

We set sample sizes as $500, 1000$ and $2000$. We take the number of replicates ($N$) as $500$. The tables contain the mean and mean Euclidean distance values of estimates. For instance, the formula for mean Euclidean distance of $\hat{\boldsymbol{\mu}}_i$ is given below

$$\|\hat{\boldsymbol{\mu}}_i - \boldsymbol{\mu}_i\| = \frac{1}{N}\left(\sum_{j=1}^{N}(\hat{\mu}_{ij} - \mu_{ij})^2\right)^{\frac{1}{2}}.$$

Similarly, the other mean Euclidean distances for the other estimates are obtained. Note that for the $\hat{\pi}_1$ the distance will be mean squared error (MSE). The formula of MSE is given

$$\widehat{MSE}(\hat{\pi}) = \frac{1}{N}\sum_{j=1}^{N}(\hat{\pi}_j - \pi)^2$$

where $\pi$ is the true parameter value, $\hat{\pi}_j$ is the estimate of $\pi$ for the $jth$ simulated data and $\bar{\pi} = \frac{1}{N}\sum_{j=1}^{N}\hat{\pi}_j$.

Table 1 shows the simulation results for the sample sizes $500, 1000$ and $2000$. In the tables, we give mean and mean Euclidean distance values of estimates and true parameter values. We can observe from this table that the proposed model is working accurately to obtain the estimates for all parameters. This can be observed from mean Euclidian distances that are getting smaller when the sample sizes increase.

**Table 1.** Mean and mean Euclidean distance values of estimates for $n = 500, 1000$ and $2000$.

| | | | Components | | | | |
|---|---|---|---|---|---|---|---|
| | | | 1 | | | 2 | |
| $n$ | Parameter | True | Mean | Distance | True | Mean | Distance |
| 500 | $\pi_1$ | 0.6 | 0.6007 | 0.0004 | | | |
| | $\mu_{i1}$ | 2 | 2.0007 | 0.2044 | -2 | -2.0186 | 0.2687 |
| | $\mu_{i2}$ | 2 | 1.9974 | | -2 | -1.9989 | |
| | $\sigma_{i,11}$ | 1.5 | 1.4928 | | 1.5 | 1.4801 | |
| | $\sigma_{i,12}$ | 0 | -0.0097 | 0.1985 | 0 | -0.0090 | 0.2561 |
| | $\sigma_{i,22}$ | 1.5 | 1.4862 | | 1.5 | 1.4835 | |
| | $\gamma_{i1}$ | 1 | 1.0030 | 0.1095 | -1 | -1.0010 | 0.1348 |
| | $\gamma_{i2}$ | 1 | 1.0053 | | -1 | -1.0029 | |
| 1000 | $\pi_1$ | 0.6 | 0.5995 | 0.0003 | | | |
| | $\mu_{i1}$ | 2 | 2.0129 | 0.1501 | -2 | -2.0028 | 0.1834 |
| | $\mu_{i2}$ | 2 | 1.9981 | | -2 | -1.9945 | |
| | $\sigma_{i,11}$ | 1.5 | 1.5004 | | 1.5 | 1.4909 | |
| | $\sigma_{i,12}$ | 0 | -0.0056 | 0.1399 | 0 | -0.0051 | 0.1740 |
| | $\sigma_{i,22}$ | 1.5 | 1.4879 | | 1.5 | 1.4942 | |
| | $\gamma_{i1}$ | 1 | 0.9987 | 0.0800 | -1 | -1.0010 | 0.0964 |
| | $\gamma_{i2}$ | 1 | 1.0022 | | -1 | -1.0024 | |
| 2000 | $\pi_1$ | 0.6 | 0.6001 | 0.0001 | | | |
| | $\mu_{i1}$ | 2 | 1.9983 | 0.1050 | -2 | -2.0022 | 0.1222 |
| | $\mu_{i2}$ | 2 | 1.9953 | | -2 | -2.0017 | |
| | $\sigma_{i,11}$ | 1.5 | 1.5011 | | 1.5 | 1.4983 | |
| | $\sigma_{i,12}$ | 0 | -0.0052 | 0.0945 | 0 | -0.0023 | 0.1155 |
| | $\sigma_{i,22}$ | 1.5 | 1.4940 | | 1.5 | 1.4940 | |
| | $\gamma_{i1}$ | 1 | 1.0004 | 0.0538 | -1 | -0.9982 | 0.0635 |
| | $\gamma_{i2}$ | 1 | 1.0026 | | -1 | -1.0011 | |

### 5.2 Real data example

In this real data example, we will investigate the bank data set which was given in Tables 1.1 and 1.2 by Flury and Riedwyl (1988) and examined by Ma and Genton (2004) to model with a skew-symmetric distribution. There are six measurements made on 100 genuine and 100 counterfeit old Swiss 1000 franc bills for this data set. This data set also analyzed by Lin (2009) to model mixtures of MSN distributions. They used the sample of $X_1$: the width of the right edge and $X_2$: the length of the image diagonal which reveals a bimodal distribution with asymmetric components. In this study, we will use the Swiss bank data to illustrate the applicability of the finite mixtures of multivariate skew Laplace distributions (FM-MSL) and also we compare the results with the finite mixtures of multivariate skew normal distributions (FM-MSN) which was given by Lin (2009). We give the estimation results in Table 2 for FM-MSN and FM-MSL. The table consists of the ML estimates, standard errors of estimates for all components, the log-likelihood, the values of the Akaike information criterion (AIC) (Akaike (1973)) and the Bayesian information criterion (BIC) (Schwarz (1978)) for the FM-MSL. Also, we give the estimation results and criterion values for FM-MSN which was computed by Lin (2009) in Table 2. According to information criterion values that the FM-MSL has better fit than the FM-MSN. In Figure 1, we display the scatter plot of the data along with the contour plots of the fitted two-component FM-MSL model. We can observe from figure that the FM-MSL captures the asymmetry and accurately fit the data.

**Table 2.** ML estimation results of the Swiss bank data set for FM-MSN and FM-MSL.

| | FM-MSN | | | | FM-MSL | | | |
| | 1 | | 2 | | 1 | | 2 | |
| | Estimate | SE | Estimate | SE | Estimate | SE | Estimate | SE |
|---|---|---|---|---|---|---|---|---|
| $w_1$ | 0.504 | 0.036 | - | - | 0.521 | 0.163 | - | - |
| $\mu_{i1}$ | 130.38 | 0.122 | 129.32 | 0.062 | 130.20 | 0.118 | 129.65 | 0.076 |
| $\mu_{i2}$ | 140.06 | 0.064 | 141.39 | 1.125 | 139.50 | 0.152 | 141.76 | 0.201 |
| $\sigma_{i,11}$ | 0.068 | 0.023 | 0.037 | 0.016 | 0.067 | 0.054 | 0.104 | 0.030 |
| $\sigma_{i,12}$ | 0.051 | 0.015 | -0.012 | 0.015 | 0.001 | 0.037 | -0.023 | 0.043 |
| $\sigma_{i,22}$ | 0.056 | 0.027 | 0.154 | 0.032 | 0.371 | 0.100 | 0.194 | 0.218 |
| $\gamma_{i1}$ | -0.230 | 0.043 | 0.494 | 0.077 | -0.017 | 0.108 | 0.034 | 0.060 |
| $\gamma_{i2}$ | -0.800 | 0.067 | 0.177 | 1.433 | 0.054 | 0.154 | -0.148 | 0.198 |
| $\ell(\hat{\Theta})$ | -310.07 | | | | **-152.30** | | | |
| AIC | 650.14 | | | | **334.60** | | | |
| BIC | 699.61 | | | | **384.08** | | | |

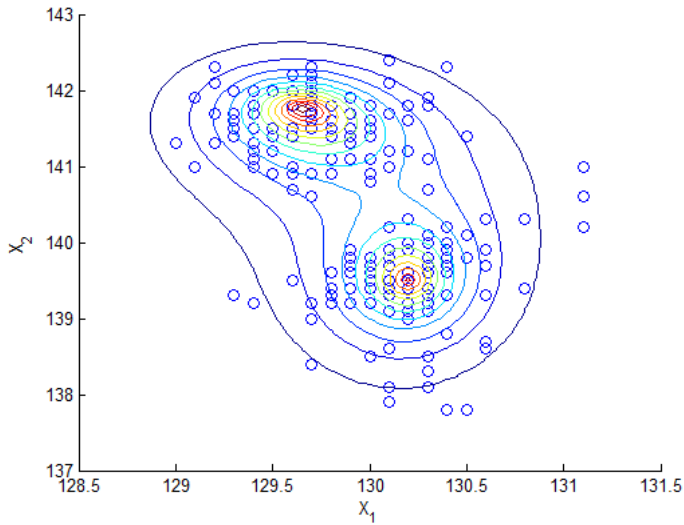

**Figure 1.** Scatter plot of the Swiss bank data along with the contour plots of the fitted two-component FM-MSL model

## 6. Conclusions

In this paper, we have proposed the mixtures of MSL distributions. We have given the EM algorithm to obtain the estimates. We have provided a small simulation study to show the performance of proposed mixture model. We have observed from simulation study that the proposed mixture model has accurately estimated the parameters. We have also given a real data example to compare the mixtures of MSL distributions with the mixtures of MSN distributions. We have observed from this example that the proposed model has the best fit according to the information criterion values. Thus, the proposed model can be used as an alternative mixture model to the mixtures of MSN distributions.

## References

Akaike, H. 1973. Information theory and an extension of the maximum likelihood principle. Proceeding of the Second International Symposium on Information Theory, B.N. Petrov and F. Caski, eds., 267-281, Akademiai Kiado, Budapest.

Arellano-Valle, R.B. and Genton, M.G. 2005. On fundamental skew distributions. Journal of Multivariate Analysis, 96(1), 93–116.
Arslan, O. 2010. An alternative multivariate skew Laplace distribution: properties and estimation. Statistical Papers, 51(4), 865–887.
Azzalini, A. and Dalla Valle, A. 1996. The multivariate skew-normal distribution. Biometrika, 83(4), 715-726.
Azzalini, A. and Capitanio, A. 2003. Distributions generated by perturbation of symmetry with emphasis on a multivariate skew t-distribution. Journal of the Royal Statistical Society: Series B (Statistical Methodology), 65(2), 367-389.
Basford, K.E., Greenway, D.R., McLachlan, G.J. and Peel, D. 1997. Standard errors of fitted means under normal mixture. Computational Statistics, 12, 1-17.
Bishop, C.M. 2006. Pattern Recognition and Machine Learning. Springer, Singapore.
Cabral, C.R.B., Lachos, V.H. and Prates, M.O. 2012. Multivariate mixture modeling using skew-normal independent distributions. Computational Statistics & Data Analysis, 56(1), 126-142.
Dempster, A.P., Laird, N.M. and Rubin, D.B. 1977. Maximum likelihood from incomplete data via the EM algorithm. Journal of the Royal Statistical Society, Series B, 39, 1-38.
Dias, J.G. and Wedel, M. 2004. An empirical comparison of EM, SEM and MCMC performance for problematic gaussian mixture likelihoods. Statistics and Computing, 14, 323-332.
Flury, B. and H. Riedwyl. 1988. Multivariate Statistics, a Practical Approach, Cambridge University Press, Cambridge.
Frühwirth-Schnatter, S. 2006. Finite Mixture and Markov Switching Models. Springer, New York.
Gupta, A. K. 2003. Multivariate skew t-distribution. Statistics: A Journal of Theoretical and Applied Statistics, 37(4), 359-363.
Gupta, A.K., González-Farías, G. and Domínguez-Molina, J.A. 2004. A multivariate skew normal distribution. Journal of multivariate analysis, 89(1), 181-190.
Hartigan, J. A. and Wong, M. A. 1979. Algorithm AS 136: A k-means clustering algorithm. Journal of the Royal Statistical Society. Series C (Applied Statistics), 28(1), 100-108.
Lin, T.I. 2009. Maximum likelihood estimation for multivariate skew normal mixture models. Journal of Multivariate Analysis, 100, 257-265.
Lin, T.I. 2010. Robust mixture modeling using multivariate skew t distributions. Statistics and Computing, 20(3), 343-356.
Lin, T.I., Ho, H. J. and Lee, C.R. 2014. Flexible mixture modelling using the multivariate skew-t-normal distribution. Statistics and Computing, 24(4), 531-546.
Ma, Y. and Genton, M.G. 2004. Flexible class of skew-symmetric distribtions, Scandinavian Journal of Statististics, 31, 459-468.
McLachlan, G.J. and Basford, K.E. 1988. Mixture Models: Inference and Application to Clustering. Marcel Dekker, New York.
McLachlan, G.J. and Peel, D. 2000. Finite Mixture Models. Wiley, New York.
Naik, D.N. and Plungpongpun, K. 2006. A Kotz-type distribution for multivariate statistical inference. In Advances in distribution theory, order statistics, and inference (pp. 111-124). Birkhäuser Boston.
Peel, D. and McLachlan, G. J. 2000. Robust mixture modelling using the t distribution. Statistics and computing, 10(4), 339-348.
Pyne, S., Hu, X., Wang, K., Rossin, E., Lin, T.I., Maier, L., Baecher-Allan, C., McLachlan, G.J., Tamayo, P., Hafler, D.A., De Jager, P.L. and Mesirov, J.P. 2009. Automated high-dimensional flow cytometric data analysis. Proc. Natl. Acad. Sci. USA 106, 8519-8524.
Sahu, S. K., Dey, D. K. and Branco, M. D. 2003. A new class of multivariate skew distributions with applications to Bayesian regression models. Canadian Journal of Statistics, 31(2), 129-150.
Schwarz, G. 1978. Estimating the dimension of a model. The Annals of Statistics, 6(2), 461-464.
Titterington, D.M., Smith, A.F.M. and Markov, U.E. 1985. Statistical Analysis of Finite Mixture Distributions. Wiley, New York.